\documentclass[11pt,letter]{article}

\usepackage{etex}
\usepackage{titlesec}
\usepackage[top=0.75in,bottom=0.75in,left=0.75in,right=0.75in]{geometry}

\titleformat{\subsubsection}[runin]{\normalfont\large\bfseries}{\thesubsubsection}{1em}{}

\usepackage{amsmath,mathrsfs,textcomp}
\usepackage{amsthm}
\usepackage{amssymb}
\usepackage{pstricks,pst-plot,psfrag}
\usepackage[all]{xy}
\usepackage{graphicx,subfigure,xspace,bm}
\usepackage{mathtools}

\usepackage{algorithm}
\usepackage[noend]{algorithmic}
\algsetup{indent=2em}

\newtheorem{theorem}{Theorem}[section]

\newtheorem{example}[theorem]{Example}

\newcommand{\real}{{\mathbb{R}}}

\usepackage{tikz}
\usetikzlibrary{external}
\usetikzlibrary{decorations.shapes,arrows,calc,decorations.markings,shapes,decorations.pathreplacing,shapes.arrows}
\usetikzlibrary{positioning}
\tikzset{main node/.style={circle,fill=blue!20,draw,inner sep=1pt},}

\newcommand{\myclearpage}{\clearpage}
\renewcommand{\myclearpage}{}

\usepackage[T1]{fontenc}
\usepackage[osf]{mathpazo}

\definecolor{BBlue}{cmyk}{.98,0.10,0,.25}

\begin{document}
                                                
\title{\rule{\textwidth}{0.4mm} {\textsc{On the lack of monotonicity of Newton-Hewer updates for Riccati equations}}\\
\rule{\textwidth}{0.4mm}}
\date{}
\maketitle
	
\vspace{-2cm}                                                

\begin{flushright}
{\footnotesize {{\bf Mohammad Akbari}}\footnote{ Department of Mathematics and Statistics at Queen's University, \texttt{13mav1@queensu.ca}. 
}\hspace{1cm}}\\
{\footnotesize {{\bf Bahman Gharesifard}}\footnote{ Department of Mathematics and Statistics at Queen's University, \texttt{bahman.gharesifard@queensu.ca}. 
}\hspace{1cm}}\\
{\footnotesize {{\bf Tamas Linder}}\footnote{ Department of Mathematics and Statistics at Queen's University, \texttt{tamas.linder@queensu.ca}. 
}\hspace{1cm}}
\end{flushright}             


\begin{abstract}%
We provide a set of counterexamples for the monotonicity of the Newton-Hewer method~\cite{GH:71} for solving the discrete-time algebraic Riccati equation in dynamic settings, drawing a contrast with the Riccati difference equation~\cite{PEC-DQM:70}. 
\end{abstract}  
	
\vspace{-0.7cm}

\section{Introduction}\label{section:intro}
This note investigates the monotonicity properties of \emph{iterative} methods for solving the discrete-time algebraic Riccati equation (DARE), which is given by
\begin{equation}\label{DARE}
P=A^\top P A-A^\top P B(B^\top P B+R)^{-1}B^\top P A + Q.
\end{equation}
As is well-known, given the discrete-time linear control system  
\[
x(t+1)=Ax(t)+Bu(t), \qquad x(0)=x_0,
\]
where $x(t)\in\real^n$ and $u(t)\in\real^m$ are the system's state and controller at time $t\geq 0$, respectively,  $A\in\real^{n\times n}, B\in\real^{n\times m}$, and the cost function 
\[
J(u)=\sum_{t=0}^\infty \Big(x(t)^\top Q x(t)+u(t)^\top R u(t)\Big),
\]
where $Q\in\real^{n\times n}, R\in\real^{m\times m}$ are positive-definite matrices, and assuming the controllability of the system $(A,B)$ and observability of $(A,Q^{1/2})$, the optimal controller which minimizes $J$ is given by 
\[
u(t)=-(B^\top P B + R)^{-1}B^\top P A,
\] 
where $P$ satisfies~\eqref{DARE}.

There are several classical iterative methods for solving the DARE in the literature, including the ones proposed in~\cite{PEC-DQM:70,GH:71}, algebraic methods~\cite{LR-PL:95}, and semi-definite programming~\cite{VB-LV:03}. In particular, these iterative methods generate a sequence of positive-definite matrices which converges to the solution of the DARE. Our main focus in this paper is on two commonly used methods, the so-called \emph{Riccati difference equation}~\cite{PEC-DQM:70}, which provably converge to the fixed point solution of the DARE, and what we call the \emph{Newton-Hewer method} which was introduced by Hewer in~\cite{GH:71}, and uses a Newton-based update to generate a sequence of positive-definite matrices which monotonically converge to the solution of the DARE when initialized at a stable policy. Let us describe these in more detail: 

The \emph{Riccati difference equation} is given by
\[
P_{t+1}=A^\top P_t A-A^\top P_t B(B^\top P_t B+R)^{-1}B^\top P_t A + Q.
\] 
It has been shown that this dynamics is monotone in terms of $P$, in the sense that $P_{t}\succeq \widehat{P}_{t}\succeq 0$ implies $P_{t+1}\succeq \widehat{P}_{t+1}\succeq 0$, see~\cite[Lemma~3.1]{CEDS:89}. Furthermore, this dynamics is monotone as a function of $(A,B,Q,R)$ in the following sense: If $P_{t}\succeq\widehat{P}_{t}\succeq 0$ and
\begin{equation}\label{eq:Hamiltonian}
\begin{pmatrix}
Q & \ A^\top \\
A & \ -BR^{-1}B^\top
\end{pmatrix}
\succeq 
\begin{pmatrix}
\widehat{Q} & \ \widehat{A}^\top \\
\widehat{A} & \ -\widehat{B}\widehat{R}^{-1}\widehat{B}^\top
\end{pmatrix},
\end{equation}
then $P_{t+1}\succeq \widehat{P}_{t+1}\succeq 0$, see~\cite{GF-GJ-HAK:96,HKW:92}. Notably and important to the discussion we will have in the next section, as long as~\eqref{eq:Hamiltonian} is satisfied, \emph{this monotonicity property holds even when the parameters $A, B, Q, R$ are time-varying}.

The \emph{Newton-Hewer method}~\cite{GH:71} is given by
\begin{align}\label{NmDare}
P_{t+1}&=A_{t}^\top P_{t+1}A_{t}+K_t^\top R K_t +Q,\nonumber\\
A_{t}&=A-B K_t,\\
K_t&=(B^\top P_t B+R)^{-1}B^\top P_t A\nonumber,
\end{align}
and it has been shown that if the system is controllable, by initializing with a stable $K_0$; i.e., $\rho(A-BK_0)<~1$, where $\rho(\cdot)$ denotes the spectral radius, $P_t$ converges monotonically, i.e., $P_1\succeq P_2 \succeq \ldots \succeq P^*$, where $P^*$ is the solution of~\eqref{DARE}.

A natural question for the Newton-Hewer dynamics is whether it has the monotonicity property that Riccati difference equation enjoys. We will show in this note that this is not the case in general, i.e., $P_{t}\succeq \widehat{P}_{t} \succeq 0$ does not necessarily imply $P_{t+1}\succeq \widehat{P}_{t+1}\succeq 0$, by providing two counterexamples, each aimed to demonstrate a facet of this lack of monotonicity.

\section{Counterexamples}
\begin{figure*}[htb!]
  \centering
    {
\includegraphics[scale=0.9]{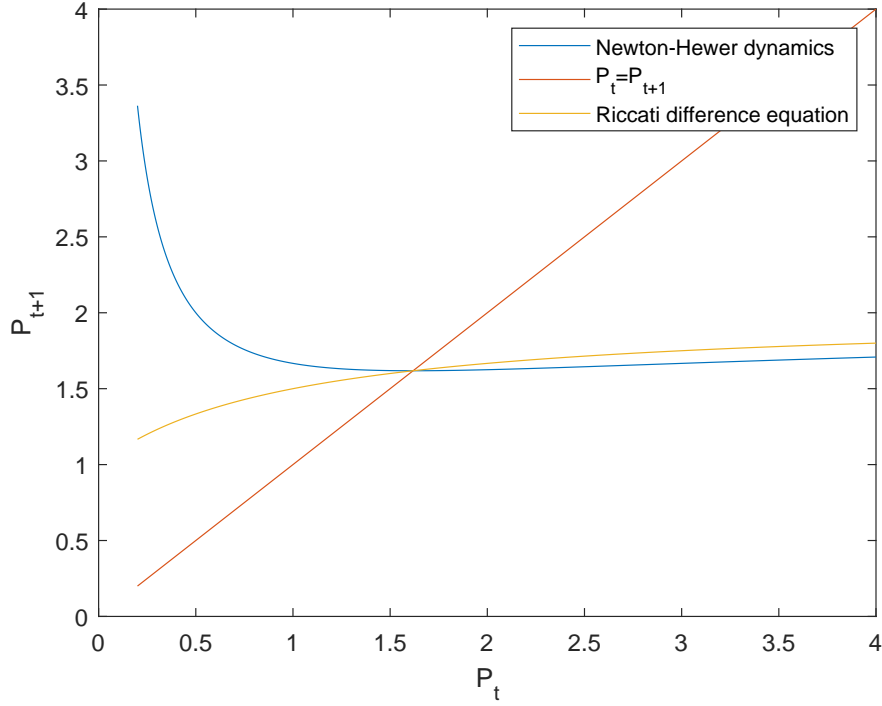}
\caption{$P_{t+1}$ as a function of $P_t$ is shown for $(A,B,R,Q)=(1,1,1,1)$ of Newton-Hewer dynamics and Riccati difference equation}
\label{Graph1}
}
\end{figure*}

The construction of our examples is done for the scalar case, and for this reason, we write the Newton-Hewer dynamics in this scenario. We assume that $ Q $ are $ R $ are positive real numbers, and do not change with time. By setting $n=m=1$, the dynamics~\eqref{NmDare} can be written as:
\[
P_{t+1}=\frac{A^2B^2P_t^2R + QB^4P_t^2 + 2QB^2P_tR + QR^2}{(P_tB^2 + R + AR)(P_tB^2 + R - AR)}.
\]
By taking derivative, it can be shown that $P_{t+1}$ as a function of $P_t$ is increasing for $P_t>P^*$ and decreasing for $P_t<P^*$, where $P^*$ is the solution to~\eqref{DARE}. For a stable policy $K_t$, $P_t$ will be larger than $P^*$ and monotonicity holds~\cite{GH:71}. We have depicted $P_{t+1}$ as a function of $P_t$ in Fig.~\ref{Graph1}, and it can be observed that the Newton-Hewer dynamics is increasing for $P_t\geq P^*$, where $P^*$ is at the intersection of the line $P_t=P_{t+1}$ and Newton-Hewer dynamics. 

We now show that if the system has time-varying $Q$ and $R$, the stabilizability properties of the controller do not necessarily imply that the system is monotone, drawing a contrast with the Riccati difference equation. 

\begin{figure*}[htb!]
  \centering
    {
\includegraphics[scale=0.9]{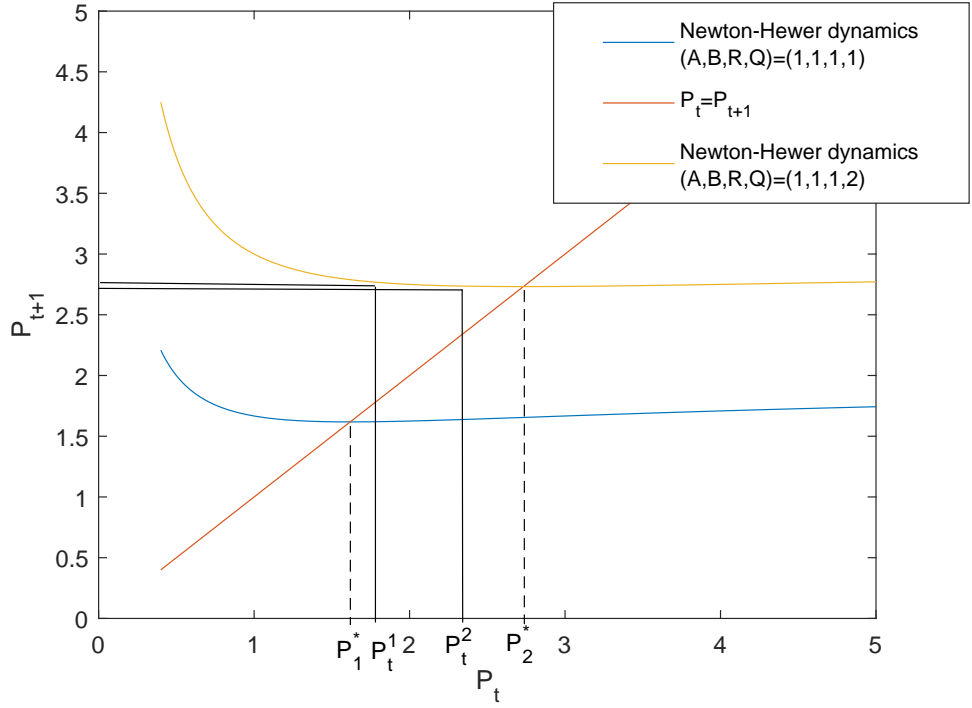}
\caption{$P_{t+1}$ as a function of $P_t$ is shown for $(A,B,R,Q)=(1,1,1,1)$ and $(A,B,R,Q)=(1,1,1,2)$ of Newton-Hewer dynamics}
\label{Graph2}
}
\end{figure*}

To this end, note that this graph depends on $A, B, R$ and $Q$, and if one of these parameters changes, $P^*$ and the graph will change. To elaborate on this, we use Fig.~\ref{Graph2} where we have depicted $P_{t+1}$ as a function of $P_t$ for two different Newton-Hewer dynamics with $(A,B,R,Q)=(1,1,1,1)$ and $(A,B,R,Q)=(1,1,1,2)$. In Fig.~\ref{Graph2}, the $P^*_1$ and $P^*_2$ refer to the solution to the DARE~\eqref{DARE} for the systems $(A,B,R,Q)=(1,1,1,1)$ and $(A,B,R,Q)=(1,1,1,2)$, respectively. If $Q_t=1, Q_{t+1}=2$ and $K_t$ are such that $P^*_1<P_t<P^*_2$, then the system for the next time step uses the orange graph 
to update $P_{t+1}$, and the reader can observe -- we prove this with carefully chosen numerical values below -- that this can lead to failure of monotonicity, i.e., $P_t\leq\widehat{P}_t$ does not necessarily imply $P_{t+1}\leq \widehat{P}_{t+1}$. Note that the system will remain monotone if $Q_{t+1}< Q_t$ for all $t$, in case the other system parameters $A, B,$ and $R$ remain fixed. 
Using this observation, we now explicitly construct the counterexample. 
\begin{example}
{\em Consider the dynamics~\eqref{NmDare}. Let the system be scalar, i.e., $n=m=1$, and let $A=1$, $B=1$, $R=1$ be fixed and $Q_t$ be time-varying. Let ${P_t}$ be the sequence generated by~\eqref{NmDare} at each time step. Given that $A,B,R$ are fixed, $P_t$ is a function of $\{Q_1,Q_2,\cdots,Q_t\}$ and $K_0$, where $K_0$ is a stable policy at time $0$. Let $\widehat{P}_t$ be the sequence generated by~\eqref{NmDare} with $A=1,B=1,R=1$ and $\widehat{Q}_t$ and $K_0$. We claim that $P_t\geq \widehat{P}_t$ does not necessarily imply that $P_{t+1}\geq \widehat{P}_{t+1}$. 

To prove this claim, we need to chose $K_0$ properly. Let 
\[
K_0=\sqrt{3}-1,
\] 
which is a stabilizing policy. Hence, by~\eqref{NmDare} we have
\[
P_1=\frac{K_0^2 R_1+Q_1}{1-(A-BK_0)^2}=\frac{4-2\sqrt{3}+Q_1}{4\sqrt{3}-6}.
\]
Given this 
\begin{align*}
K_1=&\frac{BP_1 A}{B^2 P_1+R_1}=\frac{4-2\sqrt{3}+Q_1}{2\sqrt{3}-2+Q_1},\\
P_2=&\frac{K_2^2 R_2+Q_2}{1-(A-BK_2)^2}\\
=&\frac{(8-4\sqrt{3})Q_1+(16-8\sqrt{3})Q_2+(Q_2+1)Q_1^2}{4(\sqrt{3}-1)Q_1+Q_1^2-68+40\sqrt{3}}\\
&+\frac{4(\sqrt{3}-1)Q_1Q_2+28-16\sqrt{3}}{4(\sqrt{3}-1)Q_1+Q_1^2-68+40\sqrt{3}}.
\end{align*}
Now let $Q_1=1$ and $Q_2=2$ then $P_1=1.6547$, and $P_2=2.7835$.
If we choose $\widehat{Q}_1=2$ and $\widehat{Q}_2=2$, then $\widehat{P}_1=2.7321$, and $\widehat{P}_2=2.7321$. 
This demonstrates that given $\widehat{P}_t\geq P_t$, it does not follow that $\widehat{P}_{t+1}\geq P_{t+1}$.
Fig.~\ref{Example1} shows the sequence $P_t$ (dotted line) and $\widehat{P}_t$ (dashed line) where $Q_1=1$ and $Q_t=2$ for $t\geq 2$ and $\widehat{Q}_t=2$. Furthermore, the sequence $\tilde{P}_t$ which is generated by the Riccati difference equation with initialization $\tilde{P}_1=P_1$ and the same parameters $\tilde{A}=A, \tilde{B}=B, \tilde{Q}_t=Q_t, \tilde{R}=R$ is shown (solid line) for six time steps. }

%
\begin{figure*}[htb!]
  \centering
    {
\includegraphics[scale=0.9]{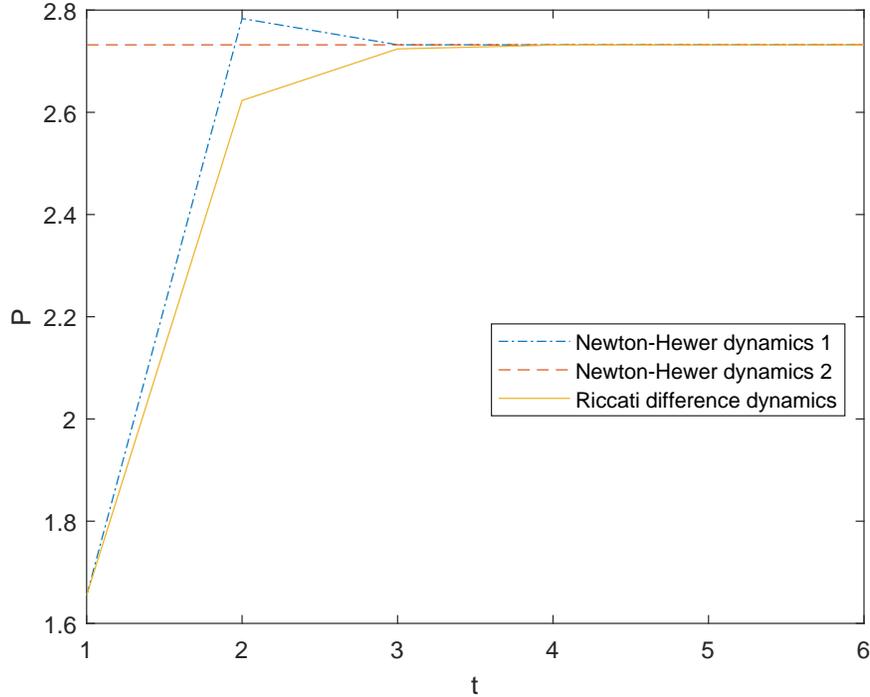}
\caption{This graph shows the Newton-Hewer dynamics for a system starts with $Q_1=1$ and $Q_t=2$ for $t>1$ (dotted line), Newton-Hewer dynamics for a system with $Q_t=2$ for all $t$ (dashed line) and  Riccati difference dynamics with $Q_1=1$ and $Q_t=2$ for $t>1$ (solid line).}
\label{Example1}
}
\end{figure*}
%
\end{example}

We conclude with providing an example which demonstrates another aspect of lack of monotonicity of Newton-Hewer dynamics. 

\begin{example}
\em{We consider two dynamics with the same $Q$ and $R$, albeit time-varying, but with different initial conditions $K_0$. 
Similar to the previous example, we assume $n=m=1$, and $A=1$, $B=1$, $R=1$ are fixed and $Q_t$ is time-varying.
We assume $Q_t$ is $1$ for odd time steps and $1.1$ for even time steps. If we choose $K_0=0.7321$ for the first system and $\widehat{K}_0=0.6180$ for the second system, we will have $P_1=1.6180$ and $\widehat{P}_1=1.6547$, and for the next time, we have $P_2=1.7351$ and $\widehat{P}_2=1.7347$, which shows that the monotonicity does not hold. Fig.~\ref{example2} illustrates the behaviour of two dynamics at the next time steps.}
\begin{figure*}[htb!]
  \centering
    {
\includegraphics[scale=0.9]{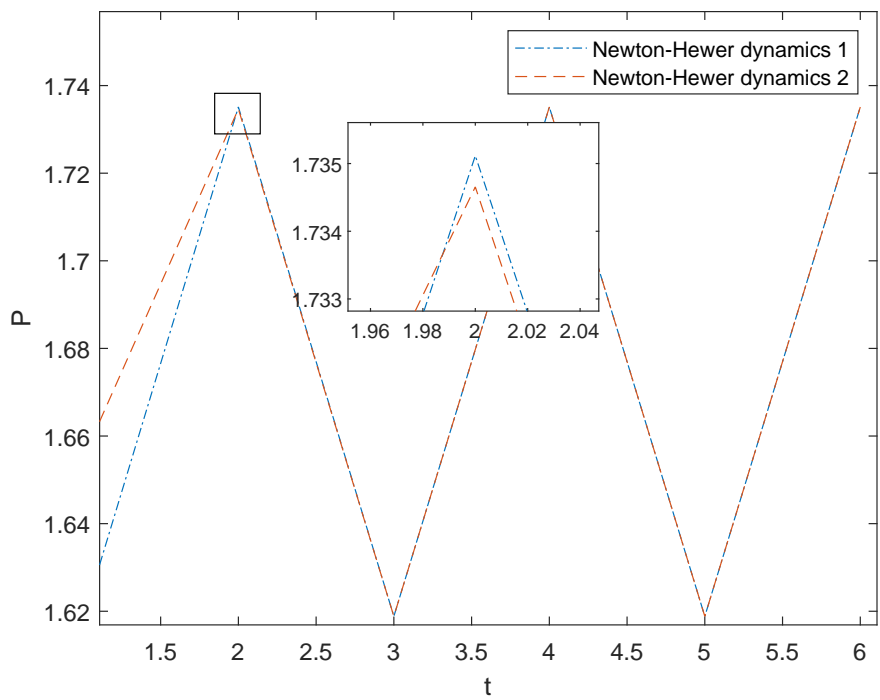}
\caption{Newton-Hewer dynamics for two systems with the same $(A,B,R,Q_t)$ and different initial condition $K_0$.}
\label{example2}
}
\end{figure*}
\end{example}

\myclearpage

\bibliographystyle{ieeetr}%

\myclearpage


\end{document}